\documentclass[11pt]{amsart}   

\usepackage{amssymb,amscd}                  

\newtheorem{thm}{Theorem}

\newtheorem{lemma}[thm]{Lemma}

\newtheorem{prop}[thm]{Proposition}

\theoremstyle{remark}

\newtheorem{remark}[thm]{Remark}

\newtheorem{example}[thm]{Example}

\newcommand{\abs}[1]{\lvert#1\rvert}

\def\U {{\mathcal U}}

\def\R {{\mathbb R}}

\def\St{{\mathrm{St}\,}}

\begin{document}


\noindent                                             
\begin{picture}(150,36)                               
\put(5,20){\tiny{Submitted to}}                       
\put(5,7){\textbf{Topology Proceedings}}              
\put(0,0){\framebox(140,34){}}                        
\put(2,2){\framebox(136,30){}}                        
\end{picture}                                         

\vspace{0.5in}

\title[Uniformities and uniformly continuous functions]
{On uniformities and uniformly continuous functions on
factor-spaces of topological groups}


\author{Siofilisi Hingano}

\address{Department of Mathematics and Statistics,
University of Ottawa, 585 King Edward Ave., Ottawa, Ontario,
Canada K1N 6N5.}


\email{shing756@science.uottawa.ca}



\subjclass[2000]{22A05}
 \keywords{ factor-space of a topological group,left uniformity, right uniformity, symmetric groups. }
\begin{abstract}
Is it true that the left and the right uniformities on a
topological group coincide as soon as every left uniformly
continuous real valued function is right uniformly continuous?
This question is known as Itzkowitz's problem, and it is still
open. We show that the generalization of the problem to
homogeneous factor-spaces of topological groups has a negative
answer.
\end{abstract}
\maketitle

The left uniform structure of a topological group $G$, which we
will denote by $\U_{l}(G),$ has a basis of entourages of the form
$$V_{l}=\{(x,y)\in G\times G: x^{-1}y\in V\},$$ where $V$ is a
neighborhood of the identity. Similarly, the right uniform
structure will be denoted by $\U_{r}(G)$ and has a basis of the
form $$V_{r}=\{(x,y)\in G\times G:xy^{-1}\in V\}.$$

A real valued function on a topological group $G$ is {\it left
uniformly continuous} if it is uniformly continuous with respect
to the left uniform structure. Similarly, for right uniform
continuity. Plainly every left uniformly continuous real valued
function is right uniformly continuous and vice versa, whenever
left and right uniform structures coincide. The converse is still
unknown. This problem is known as Itzkowitz's problem.

The question has been answered in the affirmative for particular
classes of topological groups, including locally compact groups,
metrisable groups, locally connected groups and some other classes
of groups. (see \cite {hansell, dez, kow, wit, uniform, mil,
prota}).

 Itzkowitz's question can be stated not only for topological
groups but for their homogeneous factor-spaces, because they also
support two natural uniform structures (although not necessarily
compatible with the factor topology), the left and the right ones,
see \cite {walter}. Do these two uniform structures in $G/H$
coincide as soon as every left uniformly continuous real valued
function on $G/H$ is right uniformly continuous, where $H$ is a
subgroup of the group $G$? The aim of this paper is to answer the
question on the negative, by constructing a counterexample. Notice
that on our example, both uniform structures on $G/H$ are in fact
compatible.

Let $H$ be a subgroup of a topological group $G$. Let $\pi $ be
the natural quotient mapping of $G$ onto $G/H$. The {\it left
uniform structure}, $\U_{l}(G/H),$ on $G/H$ is the finest
uniformity such that the quotient mapping
$\pi:(G,\U_{l}(G))\mapsto G/H$ is uniformly continuous. Similarly,
the {\it right uniform structure}, $\U_{r}(G/H),$ on $G/H$ is the
finest uniformity such that the mapping $\pi :(G,\U_{r}(G))\mapsto
G/H$ is uniformly continuous. It should be emphasized that if
$G/H$ is equipped with the indiscrete uniformity, then $\pi $ is
uniformly continuous. Thus the definition always makes sense.

The uniform structure $\U_{r}(G/H)$ is also called the {\em
standard uniformity} in \cite {atkin}. It is worth noticing that
the images of $\U_{r}(G)$ and $\U_{l}(G)$ under $\pi \times \pi $
are included in $\U_{l}(G/H)$ and $\U_{l}(G/H)$ respectively.

\begin{prop} \label {maile}
Let $H$ be an open subgroup of a topological group $G$. The
uniform structure $\U_{l}(G/H)$ is discrete.
\end{prop}
\begin{proof}
Since $H$ is a neighborhood  of the identity of $G$, the image of
the entourage of the diagonal $H_{\l}$ under $(\pi \times \pi )$
is the set $\{(fH,gH):f,g\in G, f^{-1}g\in H\}=\{(fH,gH):f,g\in G,
g\in fH\}=\{(fH,fH):f\in G\}$. The latter set is the diagonal of
$G/H \times G/H .$ Hence $\U_{l}(G/H)$ is just the discrete
uniform structure.
\end{proof}

 Let $X$ be a set. Let $S_{X}$ denote the symmetric group on $X$ and let $\gamma $ be a
partition of $X$. Define
$$\St_{\gamma}=\{f\in S_{X}:\forall B\in \gamma , f(B)=B\}.$$
Plainly $\St_{\gamma}$ is a subgroup of $S_{X}$.
\begin{prop} \label {losa}
The subgroups $\St_{\gamma},$ as $\gamma$ runs over all partitions
of $X$ with $\abs \gamma \leq \mathfrak c, $ form a neighborhood
basis for a Hausdorff group topology on $S_{X},$ which we will
denote by $\tau _{\frak c}.$
\end{prop}

\begin{proof}

Let $\Gamma =\{\gamma : \gamma   \text { is a partition of } $X$
\text { with }\abs \gamma\leq \mathfrak c\}.$ Since, for any
$\gamma \in \Gamma $, $\St_{\gamma}$ is a group,
$\St_{\gamma}^{2}=\St_{\gamma}=\St_{\gamma}^{-1}$.

Let $\gamma \in \Gamma $ and $g\in \St_{X}$. Define a cover $\beta
$ of $X$ as follows: $B\in \beta $ if and only if $B=g(A) $, for
some $A$ in $\gamma $. Evidently $\beta $ is in $\Gamma $ and
$g^{-1}\St_{\beta}g \subseteq \St_{\gamma}.$

If $\gamma , \beta \in \Gamma $, then define a cover $\alpha $ of
$X$ as follows: A set $D$ is in $\alpha $ if and only if $D=A\cap
B $ for some $A\in \gamma $ and $B\in \beta $. Clearly $\alpha $
is in $\Gamma $ and $\St_{\alpha}\subseteq \St_{\gamma}\cap
\St_{\beta}. $

Now the topology can be defined by taking the set
$$ \Omega =\{\St_{\gamma}: \gamma  \text { is a partition of } X,
\abs \gamma \leq \mathfrak c\},$$ as a basis of open neighborhoods
of the identity.

If $\iota $ is the identity of $S_{X} $ and $f\in S_{X}$ not equal
to $\iota $, then there exists $x \in X $ such that $f(x)\neq x.$
Put $\gamma =\{\{x\},\{f(x)\},X\setminus \{x,f(x)\}\}.$ Then
clearly $\gamma \in \Gamma $ and $f \not \in \St_{\gamma}$. This
implies $\bigcap _{\gamma \in \Gamma }\St_{\gamma}=\iota .$ Hence
$S_{X}$ is Hausdorff.

\end{proof}

\begin{remark}

\begin{enumerate}
\item The $\tau _{\frak c}$ topology on $S_{X}$ is finer than the pointwise convergence topology. For any
finite subset $M$ of $X$, let $\gamma =\{\{x\}:x\in M\}\cup
\{X\setminus M \}.$ Clearly $\gamma \in \Gamma $ and
$\St_{M}=\St_{\gamma }.$

\item If $\abs X \leq \mathfrak c$, then one obtain just the discrete topology on $S_{X}$
Indeed $\gamma =\{x\} _{x\in X}$ is a partition of $X$, of
cardinality $\leq \mathfrak c$ so $\St_{\gamma}$ is just the
identity of $S_{X}.$
\end{enumerate}

\end{remark}

Let $\gamma $ be any partition of a set $X$, with $\abs \gamma
\leq \mathfrak c $. Let $V_{\gamma}=\bigcup _{A\in \gamma}A\times
A.$ It is easy to see that the collections of all sets of the
$V_{\gamma}$ is a basis for a uniform structure on $X$. Let denote
this uniform structure by $\U _{\mathfrak c}(X)$.

\begin{prop} \label {va}
If a set $X$ and the real line $\R $ are equipped with the
uniformity $\U _{\mathfrak c} (X)$ and the additive uniformity,
respectively, then every real valued function $f$ on $X$ is
uniformly continuous.

\end{prop}

\begin{proof}

Given a function $f:X\mapsto \R, $ define a partition $\gamma $ of
$X$ as the collection of all sets $f^{-1}(x),  x\in \R .$ If now
$x,y \in X$ and $(x,y)\in V_{\gamma}$ then $f(x)=f(y)$. Hence $f$
is uniformly continuous.

\end{proof}

\begin{lemma} \label {tolu}
Let $a$ be arbitrary but fixed element of a set $X$ and let
$\gamma $ be a partition of $X$. If $b, c\in X$ then the following
are equivalent.
\begin{enumerate}
\item There are $g, f \in S_{X}$ such that $f^{-1}(\gamma
)=g^{-1}(\gamma )$ and $f(a)=b, g(a)=c.$
\item There exists $A\in \gamma $ such $b, c \in A.$
\end{enumerate}

\end{lemma}

\begin{proof}We will only prove (2) $\Rightarrow $ (1) as the other
implication is obvious. Let $f\in S_{X} $ such $f(a)=b$. Define
$g:X\mapsto X$ by
\[g(x)=\begin{cases}
c, & \mbox{if $x=a$,}\\
b, & \mbox{if $x=f^{-1}(c),$}\\
f(x), & \mbox{otherwise.}\\
\end{cases}\]
Clearly $g\in S_{X}.$ Notice that $f^{-1}(c)=g^{-1}(b),
f^{-1}(b)=g^{-1}(c) $ and for each $x\in X\setminus \{b,c\},
f^{-1}(x)=g^{-1}(x)$. It follows that that for each $A\in \gamma ,
f^{-1}(A)=g^{-1}(A)$. That is $f^{-1}(\gamma )=g^{-1}(\gamma ).$
\end{proof}

For an $a\in X$, let $St_{a}$ denote the subgroup of $S_{X}$
consisting of elements of $S_{X}$ that stabilize $a$. In our
notation this is equal to $St_{\gamma},$ where $\gamma =\{a\}\cup
\{X\setminus \{a\}\}$. Notice that every such subgroup is open in
the topology of pointwise convergence, and therefore in the
topology $\tau _{\frak c }$ as well.

\begin{thm} \label {pahulu}

Let $a\in X$ be arbitrary. Denote $H=St_{a}$. If $S_{X}$ is
endowed with the topology defined in $\tau _{\frak c}$, then
$$\U_{r}(S_{X}/H) =\U_{\mathfrak c}(S_{X}/H)\cong \U_{\mathfrak c}(X).$$

\end{thm}

\begin{proof}

The group $S_{X}$ acts on both $X$ and $S_{X}/St_{a}$ by
$(f,x)\mapsto f(x)$ and $(f,gSt_{a})\mapsto fgSt_{a},$
respectively. Define a map $\Phi :S_{X}/St_{a}\mapsto X$ by $\Phi
(fSt_{a})=f(a)$. This map can be easily shown to be well defined
and the image of $St_{a}$ under $\Phi $ is $a$.

The map $\Phi $ is equivariant, because if $g\in S_{X}$ and
$fSt_{a}\in S_{X}/St_{a}$ then
\[\Phi(gfSt_{a})=gf(a)=g(f(a))=g(\Phi (fSt_{a})).\]
 Moreover, $\Phi $ is a bijection: for every $x\in X$ there is a $f\in S_{X}$ such that
$f(a)=x$, together with the fact that every $h \in St_{a}$
stabilizes the point $a,$ imply $\Phi (fSt_{a})=x$. Also if $\Phi
(fSt_{a})=\Phi (gSt_{a})$, for some $f,g\in S_{X}$, then
$f(a)=g(a)$ so that $g^{-1}f(a)=a$. Therefore $g^{-1}f\in St_{a}$
and hence $fSt_{a}=gSt_{a}.$

Now for any partition $\gamma $ of $X$ with $\abs \gamma \leq
\mathfrak c$, the set
\begin{eqnarray*}
 \{(fH,gH):fg^{-1}\in\St_{\gamma}\}=\{(fH,gH):f^{-1}(\gamma )=g^{-1}(\gamma )\}\\ \cong
\{(f(a),g(a)):f^{-1}(\gamma )=g^{-1}(\gamma )\}.
\end{eqnarray*}
 By Lemma \ref {tolu}, the above set is equal to the set \[\{(x,y):(\exists A\in \gamma) x,y\in
A\}=\bigcup _{A\in \gamma}A\times A=V_{\gamma}.\]
\end{proof}
\begin{example}
The group $S_{X}$ on a set $X$ with cardinality greater than
$\mathfrak c$ equipped with the topology in $\tau_ {\frak c}$ with
its subgroup $H=\St_{a}$ provides a negative answer to the
generalization of Itzkowitz's question to factor-spaces.

Namely, every right uniformly continuous function on the
factor-space $S_{X}/St_{a}$ is left uniformly continuous, yet the
right and the left uniformities on $S_{X}/St_{a}$ are different.
Indeed, Proposition \ref {maile} implies that any real valued
function on $S_{X}/H$ is uniformly continuous with respect to
$\U_{l}(S_{X}/H)$. At the same time Proposition \ref {va} and
Theorem \ref {pahulu} assure that any such function is uniformly
continuous with respect to $\U_{r}(S_{X}/H)$. But Proposition \ref
{maile} and Theorem \ref {pahulu} show that these uniform
structures are different.
\end{example}
\begin{remark}
Notice that both the left and the right uniformities on $S_{X}/H$
are compatible (they generate the discrete topology on $S_{X}/H$,
which is the quotient topology). This follows from Proposition
\ref {maile} and Theorem \ref {pahulu}.
\end{remark}
\textbf{Acknowledgement}. Special thanks to my Ph.D Supervisor
Vladimir Pestov, and also Christopher J.Atkin and Peter Nickolas
for their constructive comments which helps for the product of
this paper.

\bibliographystyle{amsplain}

\end{document}